\documentclass[12pt]{amsart}
\usepackage{latexsym, amsmath,amssymb}
\usepackage{amsmath, amsthm, amssymb}
\usepackage{hyperref}
\usepackage{xcolor}
\hypersetup{
    colorlinks,
    linkcolor={red!50!black},
    citecolor={blue!50!black},
    urlcolor={blue!80!black}
}

\usepackage{graphicx}
\usepackage{float}

\providecommand{\norm}[1]{\lVert#1\rVert} 
 \numberwithin{equation}{section}

\def\XXint#1#2#3{{\setbox0=\hbox{$#1{#2#3}{%
\int}$ }
\vcenter{\hbox{$#2#3$ }}\kern-.6\wd0}}

\renewcommand{\epsilon}{\varepsilon}
\newtheorem{theorem}{Theorem}

\newtheorem{lemma}[theorem]{Lemma}
\newtheorem{corr}[theorem]{Corollary}

\newtheorem{proposition}[theorem]{Proposition}
\newtheorem{deff}[theorem]{Definition}

\newcommand{\bth}{\begin{theorem}}
\newcommand{\ble}{\begin{lemma}}
\newcommand{\bcor}{\begin{corr}}

\newcommand{\bdeff}{\begin{deff}}

\newcommand{\bprop}{\begin{proposition}}
\newcommand{\ele}{\end{lemma}}
\newcommand{\ecor}{\end{corr}}
\newcommand{\edeff}{\end{deff}}

\numberwithin{theorem}{section}

\newcommand{\eprop}{\end{proposition}}

\renewcommand{\Pi}{\varPi}

\renewcommand{\epsilon}{\varepsilon}

\begin{document}

\title[Global well-posedness of critical SQG equation on the sphere]
{Global well-posedness of critical surface quasigeostrophic equation on the sphere}
\author[D. Alonso-Or\'an]{Diego Alonso-Or\'an}
\address{Instituto de Ciencias Matem\'aticas CSIC-UAM-UC3M-UCM -- Departamento de Matem\'aticas (Universidad Aut\'onoma de Madrid), 28049 Madrid, Spain} 
\email{diego.alonso@icmat.es}
\author[A. C\'ordoba]{Antonio C\'ordoba}
\address{Instituto de Ciencias Matem\'aticas CSIC-UAM-UC3M-UCM -- Departamento de Matem\'aticas (Universidad Aut\'onoma de Madrid), 28049 Madrid, Spain} 
\email{antonio.cordoba@uam.es}
\author[A. D. Mart\'inez]{\'Angel D. Mart\'inez}
\address{Instituto de Ciencias Matem\'aticas (CSIC-UAM-UC3M-UCM) -- Departamento de Matem\'aticas (Universidad Aut\'onoma de Madrid), 28049 Madrid, Spain} 
\email{angel.martinez@icmat.es}

\begin{abstract}
In this paper we prove global well-posedness of the critical surface quasigeostrophic equation on the two dimensional sphere building on some earlier work of the authors. The proof relies on an improving of the previously known pointwise inequality for fractional laplacians as in the work of Constantin and Vicol for the euclidean setting.
\end{abstract}

\maketitle

\section{\textbf{Introduction}}\label{s:intro}


In this paper we prove global well-posedness of the critical surface quasigeostrophic equation on the two dimensional sphere. The proof relies on an integral representation of the fractional Laplace-Beltrami operator on a general compact manifold and an instantaneous continuity result for weak solutions (cf. Theorems \ref{1} and \ref{prop} stated below). The representation has a benign error term which allows an improvement of the C\'ordoba-C\'ordoba inequality \cite{CC2} in the spirit of Constantin and Vicol (cf. \cite{AOCM3, CVi}). Their work handles global existence of the critical surface quasigeostrophic equation in $\mathbb{R}^n$, which followed landmark results obtained independently by Kiselev, Nazarov and Volberg \cite{KNV} and Caffarelli and Vasseur \cite{CV}. Underneath our arguments below we exploit the rich group of isometries of the sphere. 

\subsection{Global well-posedness of the critical surface quasigeostrophic equation on the sphere}\label{sqg}

The critical surface quasigeostrophic equation on the sphere has been treated in \cite{AOCM} and is given by
\begin{equation}\label{eq:SQG}
\begin{cases}	
\theta_{t}+u \cdot \nabla_{g} \theta+\Lambda\theta=0 \\
u = \nabla_{g}^\perp \Lambda^{-1} \theta\\
\theta(0)=\theta_0
\end{cases}
\end{equation}
One studies the evolution of some class of initial data. In their previous work \cite{AOCM} the authors established an explicit modulus of continuity for weak solutions:

\begin{theorem}\label{1}
Given an initial datum $\theta_0\in L^2(\mathbb{S}^2)$ any weak solution of (\ref{eq:SQG}) becomes instantanously continuous with an explicit modulus of continuity, for any time $t>0$.
\end{theorem}

More especifically, if $t\geq t_{0}>0$ then the modulus of continuity is shown to have the form $\omega(\rho)=O\left((\log(1/\rho))^{-\alpha}\right)$ for some fixed $\alpha=\alpha(t_0)>0$ which degenerates as $t_0$ tends to zero. The proof exploits Caffarelli and Vasseur's analysis which is based on De Giorgi's techniques on each scale at a time (cf. \cite{AOCM} for details). In the present paper we provide also some regularity results for this equation

\begin{theorem}\label{global}
There is global well posedness in $H^s(\mathbb{S}^{2})$ for any $s> 3/2$. In fact, any solution with such initial datum becomes smooth instantanously.
\end{theorem}

An inmediate consequence is the following result (analogous to that of Nazarov, Kiselev and Volberg \cite{KNV}).

\begin{theorem}
Given an initial data $\theta_0\in C^{\infty}(\mathbb{S}^2)$, the solution will remain smooth for all times $t>0$.
\end{theorem}

The proof will follow closely the strategy of Constantin and Vicol \cite{CVi} which is based on the nonlinear maximum principle established with help of the aforementioned explicit integral representation for the fractional Laplace-Beltrami operator. The main difference with their exposition is that we know a priori that the solution is continuous uniformly for times $t\geq t_0>0$, invoking Theorem \ref{1}. This result is of independent interest and allows us to let aside their technical stability result, which might nevertheless be true in $\mathbb{S}^n$. Observe that the modulus of continuity is not uniform for small times. As a consequence, we need to close a gap for small times using local existence which follows by a standard energy estimate (cf. Appendix \ref{s:appendix}). Theorem \ref{1} implies what Constantin and Vicol define as {\em only small shocks} condition. This condition allows to control certain contributions at mesoscopic scales.
uniformly in time $t\geq t_0>0$. This kind of control is needed in the proof of the following result:

\begin{proposition}\label{gradientcontrol}
 Let $\theta(x,t)$ be a weak solution of the critical surface quasi-geostrophic equation (\ref{eq:SQG}) then it satisfies the following bound
 \[ \sup_{t\geq t_0}\norm{\nabla_g \theta}_{L^{\infty}(\mathbb{S}^{2})} \leq C\left(\norm{\theta_{0}}_{L^{\infty}(\mathbb{S}^{2})},\norm{\nabla_g\theta_{0}}_{L^{\infty}(\mathbb{S}^{2})}\right)\]
\end{proposition}
 
Notice this gradient control inmediately implies Theorem \ref{global} (cf. Appendix \ref{s:appendix}). It is precisely at this point where we need to be restricted to the two dimensional sphere. However, similar arguments deal with arbitrary dimension tori or the euclidean space, as has already been treated in the literature (cf. \cite{CV}, \cite{CVi}, \cite{KNV}).

\subsection{The integral representation}

As usual let $(M,g)$ be a compact manifold of dimension $n\geq 2$ whose Laplace-Beltrami operator is denoted by $-\Delta_g$. Let us state the representation formula now

\begin{theorem}\label{prop}
Let $f$ be smooth and $s\in(0,1)$, for any fixed $N$ and diagonal cut off function $\chi$ one has the following
\[(-\Delta_g)^{s}f(x)=\textrm{P.V.}\int_M\frac{f(x)-f(y)}{d(x,y)^{n+2s}}(c_{n,s}\chi u_0+k_N)(x,y)d\textrm{vol}_g(y)+O(\|f\|_{H^{-N}(M)})\]
where $k_N(x,y)=O(d(x,y))$ is some smooth function, the implicit constant depends on $N$ and $c_{s,n}>0$ is a constant independent of $N$ and $u_0$ is a smooth function with $u_0(x,x)=1$.
\end{theorem}

Notice that the norm in the error might be taken to be $L^{\infty}$ for $N$ big enough. For further details of its proof we refer the reader to \cite{AOCM3}. Let us remark that the explicit smoothing effect of the error term is crucial for the application this paper deals with.  One may compare this with similar expresions in the case of flat tori or euclidean space (cf. \cite{CC,CC2,S1}). We carry the proofs for arbitrary $n$ whenever it does not affect the exposition clearness.

\section{\textbf{Pointwise estimates}}\label{s:point}

In this section we provide several observations that will be instrumental in the sequel. In particular we prove pointwise commutator estimates related to the ones appearing in the work of Constantin and Ignatova (cf. \cite{CI1,CI2}) which allow to estimate $[\Lambda,\nabla]$. This commutator involves the action of $\Lambda$, a pseudodifferential operator, on fiber bundles, which can be defined in several ways. However, our pointwise estimates do not seem to be an inmediate consequence of this general setting.


In the following some needed explicit calculations will be presented, we shall reduce our exposition of them to the two dimensional sphere although the arguments work in arbitrary dimension. Let us denote by $R_s$ a rigid rotation of the sphere around some axis where $s$ denotes the arc length of the particles moving in the corresponding equator and by $\dot{R}_{s}$ the vector field it generates infinitesimaly. Since rotations are isometries, their generator commutes with the laplacian which easily implies that $[\Lambda,\dot{R}_s]$ vanishes.

\begin{lemma}\label{conmutador2}
Let $\alpha\in (0,2)$ and suppose that the smooth function $a$ satisfies $a(x)-a(y)=O(d(x,y)^2)$. Then there is a constant $C=C(n,\alpha)$ such that the following pointwise commutator estimate holds
\[ |[\Lambda^{\alpha}, a]f(x)|\leq C\|f\|_{\infty}\]
\end{lemma}

\textsc{Proof of Lemma \ref{conmutador2}:} Let us employ the following representation of the fractional operator (up to a constant)
\[\Lambda^{\alpha} f(x)=\int_0^{\infty}t^{-1-\alpha/2}\left(f(x)-\int_MG(x,y,t)f(y)d\textrm{vol}_g(y)\right)dt\]
from which it easily follows that the commutator satisfies
\[[\Lambda^{\alpha},a]f(x)=\int_0^{\infty}t^{-1-\alpha/2}\int_M G(x,y,t)(a(y)-a(x))f(y)d\textrm{vol}_g(y)dt\]
Which might be estimated for small times quite crudely employing the following upper bound for the heat kernel (cf. \cite{LY}, Corollary 3.1)
\[G(x,y,t)\leq C(M,g)\frac{e^{-\frac{d_g(x,y)^2}{5t}}}{t^{n/2}}\]
The proof concludes observing that one can also estimate the rest easily taking advantage of the exponential decay of the heat kernel on compact manifolds for large times.

Alternatively one may actually provide a proof just computing the commutator using our kernel representation (Theorem \ref{prop}). 

Before proceeding any further let us include the following general result which will be employed in the sequel.

\begin{lemma}\label{lem:positive}
Let $f\geq 0$ be some smooth function on a compact riemannian manifold $(M,g)$ and denote by $x\in M$ the point where it reaches it maximum. Then, for any $\alpha\in(0,2)$
\[\Lambda^{\alpha} f(x)\geq 0\]
\end{lemma}

This is somehow suprising since curvature might have some effects in view of our representation formula. It is nevertheless true in the stated generality as we will prove now (cf. \cite{CM}).

\textsc{Proof of Lemma \ref{lem:positive}:} let us introduce the following Cauchy problem for a fractional heat equation on the manifold, namely
\[\left\{\begin{array}{l}
\frac{d}{dt}h=-\Lambda^{\alpha}h\\
h(\cdot, 0)=f
\end{array}\right.\]
It is well known that it satisfies the following maximum principle \[\|h(\cdot,0)\|_{L^{\infty}_x}\geq\|h(\cdot, t)\|_{L^{\infty}_x}.\] As a consequence  at that maximum point $x$ we have $h(x,t)-h(x,0)\leq 0$, dividing then by $t$ and letting $t$ approach zero one gets $\frac{d}{dt}h(x,0)\leq 0$. But this is equivalent to our claim.

We include now an approximation of which we will take advantage in the next section and although it holds in general dimensions, for the sake of simplicity we shall present the details of the proof only in dimension two. We will show that for any point $x$ we can approximate to second order the infinitesimal rotations $\dot{R}_1$, $\dot{R}_2$ corresponding to the rotations induced by a given orthonormal system of vectors in $T_xM$ with vector fields $\partial_1$ and $\partial_2$ in some appropiate coordinates. 

Let us consider the stereographic projection with $p$, its south pole, the origin of coordinates $(0,0,0)\in\mathbb{R}^3$. 

\[(w_1,w_2)=\left(\frac{x}{1-z},\frac{y}{1-z}\right)\]

\begin{figure}[htb]
\centering
\includegraphics[width=90mm]{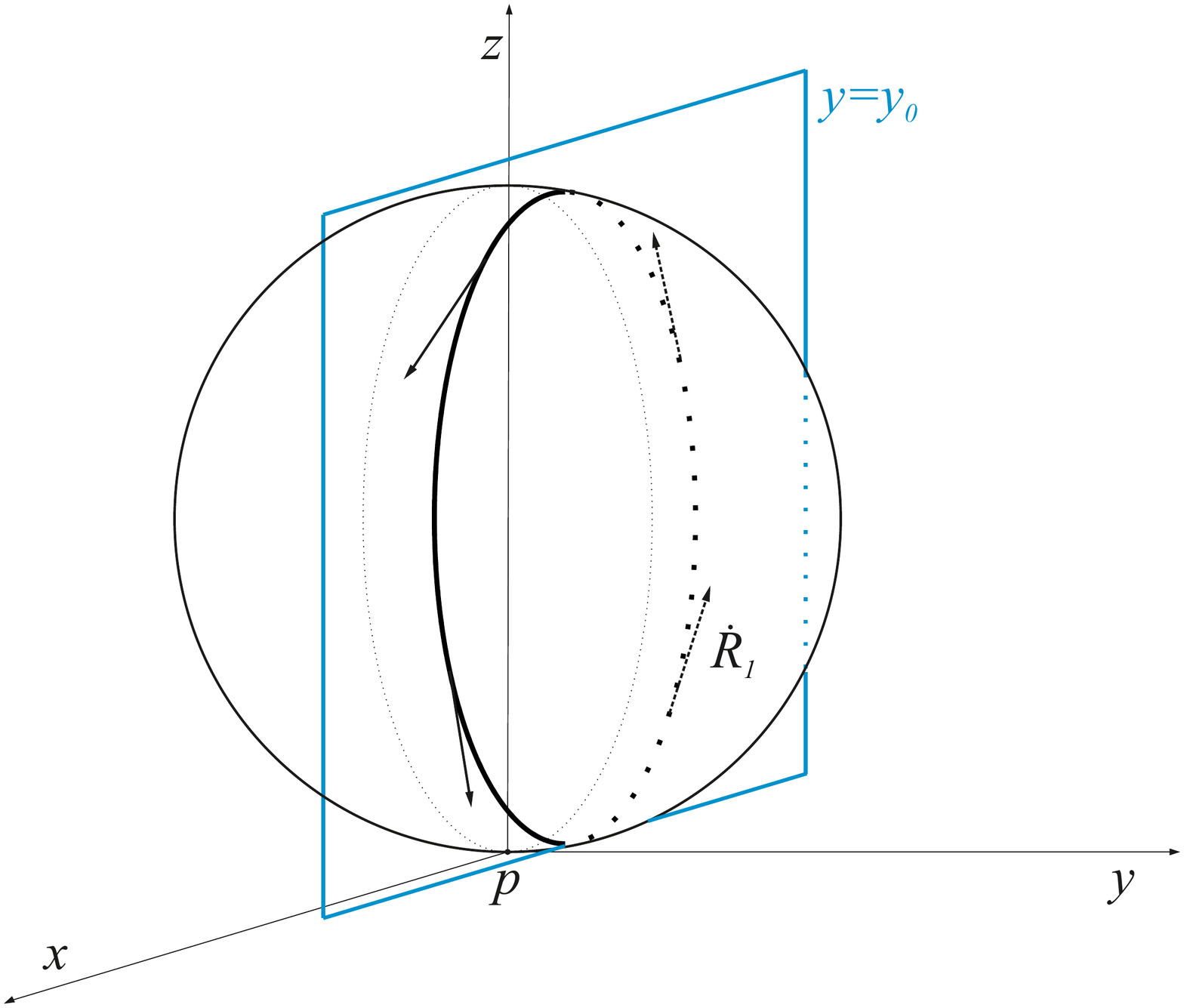}
\caption{Rotations on the sphere}\label{fig:1}
\end{figure}

Next, we compute $\dot{R}_1$ in this coordinates. For that purpose we parametrize the circle of rotation of some point near the south pole, corresponding to $y=y_0$ constant (see figure \ref{fig:1}).

Namely
\[\left(\sqrt{1-y_0^2}\sin(\alpha),y_0,\sqrt{1-y_0^2}\cos(\alpha)-1\right)\]
then $\dot{R}_1$ corresponds to derivative with respect to $\alpha$, which, in stereographic coordinates corresponds to
\[\left(\frac{\sqrt{1-y_0^2}\sin(\alpha)}{2-\sqrt{1-y_0^2}\cos(\alpha)},\frac{y_0}{2-\sqrt{1-y_0^2}\cos(\alpha)}\right).\]
Straighforward differentiation and use of Taylor approximations of the functions therein implies
\[\dot{R}_1=(1+O(h^2))\partial_{w_1}+O(h^2)\partial_{w_2}\]
for any $y_0,\alpha\leq h$. Let us denote by $a^{ij}$ the coefficients of the change of coordinates. Similarly for $\dot{R}_{2}$. As a consequence $\dot{R}_{i}$ can coincides with $\partial_i$ up to an error of second order. Finally, one have to compute the metric tensor in this coordinates, $g_{ij}$, and the same method shows that it is a perturbation of second order of the identity, i.e. $g_{ij}(y)=\delta_{ij}(x)+O(|x-y|^2)$). This fact will be very convenient in order to apply our previous Lemma ~\ref{conmutador2} effectively.

One may observe also that the stereographical projection coordinates are not far from being equal to the polar coordinates. That is, they differ only on a second order perturbation, allowing to transfer many estimates and properties from one to the other.
\section{\textbf{Nonlinear lower bounds}}


Here we prove in our setting a refinement of the pointwise inequality obtained by Constantin and Vicol. Their work happens in euclidean space $\mathbb{R}^n$ which is an unbounded domain, however since we are interested in compact manifold some extra hypothesis have to be imposed in the statement together with the curvature effects in the error term. From now on we will employ the notation $\Lambda=(-\Delta_g)^{1/2}$.
 	
\begin{proposition}\label{muesca}
Let $f$ be a smooth function on the sphere and $0<\alpha<2$. Then, provided $|\nabla_g f(x)|\geq C\|f\|_{\infty}$, we have the pointwise bound 
\[ \nabla_g f(x)\cdot\nabla_g\Lambda^{\alpha} f(x) \geq \frac{1}{2} \Lambda^{\alpha}(|\nabla_g f|^{2})(x)+ \frac{1}{4}D(x)+\frac{|\nabla_g f(x)|^{2+\alpha}}{c\|f\|_{\infty}^{\alpha}}+O(\|\nabla_gf\|_{\infty}^2) \]
where $D$ denotes some positive functional to be defined in the proof and the positive constants depend on $n$ or $\alpha$ but it is independent of $x$. 
\end{proposition}

The idea of the proof is to employ the representation formula mentioned above and use the commutator (Lemma \ref{conmutador2}) in order to obtain positivity in the principal term. We do take advantage of the natural isometries of the round spheres. Direct use of microlocal analysis does not seem to help much at this point.

\textsc{Proof of Proposition \ref{muesca}:} fix the point $x$ and around it consider the stereographical coordinates introduced in the previous section. Then we have $g^{ij}(x)=\delta^{ij}(x)$ and using our previous section one may compute the left hand side of the inequality as follows:
\[\begin{split}
g^{ij}(x)\partial_if(x)g_{jk}(x)g^{k\ell}(x)\partial_{\ell}\Lambda^{\alpha} f(x)&=\partial_i f(x)g^{i\ell}(x)\partial_{\ell}\Lambda^{\alpha} f(x)\\
&=\dot{R}_if(x)a^{i\ell}(x)\dot{R}_{\ell}\Lambda^{\alpha} f(x)\\
&=\dot{R}_if(x)a^{i\ell}(x)(\Lambda^{\alpha} \dot{R}_{\ell}f)(x).\\
\end{split}\]
That is:
\[\nabla_g f(x)\cdot\nabla_g\Lambda^{\alpha} f(x)=\nabla f(x)\cdot\nabla\Lambda^{\alpha}f(x)\]
since the distorsion $a^{ij}$, which is is a perturbation of second order of the metric, is not noticed at $x$. Let us introduce for simplicity the following notation $\nabla f=(\dot{R}_1f,\ldots,\dot{R}_nf)$. The representation formula given in Proposition \ref{prop} for the fractional laplacian yields
\[ \nabla f(x)\cdot \Lambda^{\alpha}\nabla f(x)=\frac{1}{2} \Lambda^{\alpha}(|\nabla f|^{2})(x) + \frac{1}{2} D(x)+ E(x)+O(\|\nabla f\|_{\infty}^2)\]
where 
\[ D(x) = c_{\alpha,n} P.V.\int_{\mathbb{S}^n} \frac{|\nabla f(x)-\nabla f(y)|^{2} }{d(x,y)^{n+\alpha}}u_0\chi(x,y)d\textrm{vol}_{g}(y) \]
and the error terms are
\[ E(x) = \frac{1}{2} \int_{\mathbb{S}^n} \frac{|\nabla f(x)-\nabla f(y)|^{2}}{d(x,y)^{n+\alpha}}k_{N}(x,y) d\textrm{vol}_{g}(y)\]
one might wish to subsume this in, say, $\frac{1}{8}D(x)$, which will be possible due to the nature of $k_{N}$. Indeed, it is clear from the construction of the parametrix that $k_N$ is supported where $\chi(x,y)$ is supported and it is bounded a priori by some constant, say, $M$ independently of $\chi$ but depending on $N$. As a consequence one may choose the cut off $\chi$ to be supported in $d_g(x,y)\leq \frac{1}{8} M$ and that choice provides the desired estimate. Therefore, we can absorb that term in $D(x)$ after changing conveniently the constant.

Using the commutator estimates from the Lemma \ref{conmutador2} and an appropiate cut off function we have
\[\Lambda^{\alpha}(|\nabla f|^2)(x)-\Lambda^{\alpha}(|\nabla_g f|^2)(x)=O(\|\nabla_g f\|^{2}_{\infty}).\]
We will be done if we prove,
\[ D(x)\geq \frac{|\nabla f(x)|^{2+\alpha}}{c\norm{f}^{\alpha}_{L^{\infty}}}.\] To do that let us introduc a smooth cut off $\eta_{\rho}$ supported outside a ball of radius $2\rho$, say, and equal to zero inside the ball of radius $\rho$ around $x$. We will optimize the radius $\rho$, provided that it is smaller than the injectivity radius, to get the desired inequality. Let us observe now that
\[ D(x) \geq c_{\alpha,n} \int_{\mathbb{S}^{n}} \frac{|\nabla f(x)-\nabla f(y)|^{2}}{d(x,y)^{n+\alpha}} \eta_{\rho}(y)\chi(x,y) d\textrm{vol}_{g}(y).\]
Next we use the fact that for every $y$ we have:
\[ |\nabla f(x)-\nabla f(y)|^{2} \geq |\nabla f(x)|^{2} - 2\nabla f(x)\cdot\nabla f(y). \]
Hence we get
\[\begin{split}
 D(x) &\geq c_{\alpha,n}|\nabla f(x)|^{2}\int_{\mathbb{S}^n} \frac{\eta_{\rho}(y)\chi(x,y)}{d(x,y)^{n+\alpha}}d\textrm{vol}_{g}y \\
 &\quad\qquad\qquad\qquad\qquad\qquad\qquad- c_{\alpha,n}|\nabla f(x)|\left|\int_{\mathbb{S}^n}\frac{\nabla f(y)\eta_{\rho}(y)\chi(x,y)}{d(x,y)^{n+\alpha}} d\textrm{vol}_{g}(y)\right| \\
      &\geq c_{\alpha,n}C|\nabla f(x)|^{2}\int_{\mathbb{S}^n} \frac{\eta_{\rho}(y)\chi(x,y)}{d(x,y)^{n+\alpha}}d\textrm{vol}_{g}y \\
     & \quad\qquad\qquad\qquad\qquad\qquad\qquad- c_{\alpha,n}|\nabla f(x)|\norm{f}_{L^{\infty}}\int_{\mathbb{S}^n}\left| \nabla\left(\frac{\eta_{\rho}(y)\chi(x,y)}{d(x,y)^{\alpha+n}}\right)\right| d\textrm{vol}_{g}(y) \\
      &\geq C_{1}\frac{|\nabla_g f(x)|^{2}}{\rho^{\alpha}}-C_{2}\frac{|\nabla_g f(x)|\norm{f}_{L^{\infty}}}{\rho^{\alpha+1}}-\frac{|\nabla_g f(x)|^2}{M^{\alpha}}
\end{split}\]
where $C_{1},C_{2}$ might depend on $\alpha$, the dimension and the cut off, which are already fixed. We would like to set \[\rho=\frac{C_{2}\norm{f}_{L^{\infty}}}{2C_{1}|\nabla_g f(x)|},\] in order to obtain the nonlinear bound. But one must take care that such a ball might be chosen first and that is only possible if
\[ |\nabla_g f(x)|\geq C_{3}\norm{f}_{L^{\infty}} \]
where $C_3$ depends on the previous constants and the diagonal support width of $\chi$. 

\section{\textbf{Proof of Proposition ~\ref{gradientcontrol}}}

We will suppose that $\theta$ is smooth in $[0,T)$, otherwise we regularize it introducing some artificial viscosity $\nu\Delta_g\theta$ term and proving that the estimates do not depend on $\nu$. Then one takes the vanishing viscosity limit $\nu\rightarrow 0$. Some comments are in order: studying the evolution of $|\nabla_g\theta|^2$ a term of the form $\nabla_g\Delta_g$ appears which is troublesome but, fortunately, there is a way to overcome this difficulty. Namely: expressing the gradient as a combination of the form $a^{ij}(x)\dot{R}_j$. The rotations commute with the Laplace-Beltrami operator and the commutator will have terms involving first derivatives of $a^{ij}(x)$, which vanish, and second derivatives of $a^{ij}(x)$, which are uniformly bounded, coupled with (local) first derivatives of $\theta$, which can be absorved by $|\nabla_g\theta|(x)$ (cf. Proposition \ref{conmutador2} above). The rest of the argument is standard.

A priori the limit functions might correspond to different weak solutions but one can prove that any two weak solutions $\theta_1$ and $\theta_2$ coincide provided that one of them is smooth (which we know a posteriori!). Another fact we will be building in is that an a priori estimate on $\norm{\nabla_g \theta(\cdot, t) }_{L^{\infty}}$ inmediatly implies that $\theta$ is smooth. This rather elementary facts are included for the sake of completeness in Appendix \ref{s:appendix}.

Applying $\nabla_{g}$ to the equation and taking the scalar product with $\nabla_{g}\theta$ yields
\[\frac{1}{2}(\partial_{t} + u\cdot \nabla_{g})|\nabla_{g}\theta|^{2} + \nabla_g\theta\cdot\nabla_g\Lambda_{g}\theta+\nabla_gu:\nabla_g\theta\cdot\nabla_g\theta+u\cdot\nabla_g\theta\cdot\nabla_g\theta=0\]
Now if we could apply the pointwise lower bound of Proposition \ref{muesca} we would get:
\[\frac{1}{2}L(|\nabla_{g}\theta|^{2}(x))+\frac{1}{4}D(x)+ \frac{|\nabla_g \theta(x)|^{3}}{c\norm{\theta}_{L^{\infty}}} \leq O(|\nabla_{g}u(x)||\nabla_{g}\theta(x)|^{2}+|u(x)||\nabla_g\theta(x)|^2)+O(\norm{\nabla_{g}\theta}^{2}_{L^{\infty}})\]
where $L(f)=\partial_{t}f+u\cdot\nabla_{g}f+\Lambda_{g} f$ and $c$ is a positive universal constant. This holds true provided $|\nabla_g\theta|(x)\geq C\|\theta\|_{L^{\infty}}$. We claim that an estimate of the form
\begin{equation}
\frac{1}{2}L(|\nabla_{g}\theta|^{2})(x)+\frac{|\nabla \theta(x)|^{3}}{c\norm{\theta}_{L^{\infty}}} \leq C |\nabla_{g}\theta(x)|^{2}+ O(\norm{\nabla_{g}\theta}^{2}_{L^{\infty}})\label{test}
\end{equation}
is valid for some constant $C=C(\epsilon,M,\norm{\theta}_{L^{\infty}})$. 

The result then follows intuitively reading this estimate at a point where $|\nabla_{g}\theta(x)|$ attained its maximum. Indeed, as a consequence of the positivity of the fractional Laplacian and that the gradient vanishes when evaluated on a maximum point, one expects that $\frac{d}{dt}|\nabla_g\theta|^2\leq0$ which obstructs its indefinite growth. To finish we observe that in the case that we can not apply the nonlinear bound \ref{muesca}, we would have $\|\nabla_g\theta\|_{\infty}\leq C\|\theta\|_{\infty}$ which is even better. We postpone the rigorous argument to the end of the section.

Next step is to show an integral representation of the main term of the velocity $u$ similar to the euclidean one (cf. \cite{CVi}). We will take advantage of the sphere geometry but we believe that it could be done in a more general setting. 

First one may express in local coordinates around $x$ the derivative
\[\begin{split}
\partial_{\ell}u_i(x)&=\partial_{\ell}(g_{\perp}^{ij}\partial_j\Lambda_{g}^{-1}\theta)(x)=a_{\ell k}(x)\dot{R}_k(g_{\perp}^{ij}a_{jn}\dot{R}_n\Lambda_{g}^{-1}\theta)(x)\\
&=a_{\ell k}(x)g_{\perp}^{ij}(x)a_{jm}(x)\Lambda_{g}^{-1}\dot{R}_k\dot{R}_m\theta(x)+a_{\ell k}(x)\dot{R}_kg_{\perp}^{ij}(x)a_{jm}(x,t)\Lambda_{g}^{-1}\dot{R}_m\theta(x)
\end{split}\]
Next, using the integral representation and integration by parts we get the following estimate:
\[|\partial_{\ell}u_i(x)|\leq C\left|\int_{\mathbb{S}^n}\frac{\dot{R}_k\dot{R}_m\theta(y)}{d(x,y)^{n-1}}(c_{n}u_{0}\chi+k_N)(x,y)d\textrm{vol}_g(y)\right|+\gamma\|\nabla_g\theta\|_{\infty}+O(\log\gamma\|\theta\|_{\infty})\]
for some constant $C$ and a parameter $\gamma$ to be chosen later to be small enough. The harmless error terms in the right hand side come from
\[O(1)\int_{\mathbb{S}^n}\frac{\dot{R}_m\theta(y)}{d(x,y)^{n-1}}(c_{n}u_{0}\chi+k_N)(x,y)d\textrm{vol}_g(y)+O(\|\theta\|_{\infty}).\]
Indeed, the integral might be splitted smoothly in two parts: near points, say $d(x,y)\leq \gamma$; points at distance $\gamma$; and far away points $d(x,y)> \gamma$. The first is bounded by $\gamma\|\nabla_g\theta\|_{\infty}$ and we might choose $\gamma$ as small as we wish so that the cubic term is absorved by the one in the left hand side of the inequality (\ref{test}). The remainder is bounded after integration by parts by $\log(\gamma)\|\theta\|_{\infty}$ which nevertheless might be a rather large constant. The $O(|u||\nabla_g\theta|^2)$ term is handled similarly. 

Let us now deal with the main part of the estimate which after integration by parts reads
\[\int_{\mathbb{S}^n}\dot{R}_m\theta(y)\dot{R}_k\left(\frac{c_{n}\chi(x,y)+k_N(x,y)}{d(x,y)^{n-1}}\right)d\textrm{vol}_g(y).\]
We split smoothly the integral into three summands, namely:
\[I_{\textrm{in}}(x)+I_{\textrm{med}}(x)+I_{\textrm{out}}(x).\]
an inner piece for near points $d(x,y)< \rho$, for some specific $\rho(x)>0$ to be choosen later, a middle piece $\rho<d(x,y)<\epsilon$ (where $\epsilon>0$ will also be choosen later) and an outer part for points $d(x,y)>\epsilon$. 

To bound the inner piece, we notice that writting it in stereographic coordinates the most singular term has the following form:
\[\int_{d(x,y)<\rho} \dot{R}_k\theta(y)a^{im}(y)\partial_m\left(\frac{h(x,y)}{|x-y|^{n-1}}\right)d\textrm{vol}_g(y)\]
for some nice function $h$  (which might be computed explicitly in this spherical case) and where we understand that there is a smooth cut off in the integral. This expression allows us to take advantage of the cancellation in spheres around $x$ to achieve
\[\int_{d(x,y)<\rho}\frac{|\dot{R}_k\theta(y)-\dot{R}_k\theta(x)|}{|x-y|^n}h(x,y)d\textrm{vol}_g(y)+O(\rho(x,y)\|\nabla_g\theta\|_{\infty}),\]
where the second term above comes from the estimations of lower order terms. Now applying Cauchy-Schwarz inequality we get 
\[ |I_{\textrm{in}}(x)| \leq C \sqrt{D(x)\rho(x)} + O(\rho(x,y)\|\nabla_g\theta\|_{\infty}).\]
Therefore choosing $\rho(x)=\frac{C\norm{\theta}_{L^{\infty}}}{|\nabla_{g}\theta(x)|}$ for suitable $C$, and using Young's inequality we obtain
\[|I_{\textrm{in}}(x)||\nabla_{g}\theta(x)|^{2}\leq \frac{1}{8}D(x) + \frac{|\nabla_{g}\theta|^{3}}{c\norm{\theta}_{L^{\infty}}}+ O\left(\norm{\nabla_{g}\theta}^{2}_{L^{\infty}} \right). \]
The medium piece can be rewriten as:
\[\int_{\epsilon>d(x,y)>\rho}\dot{R}_m(\theta(y)-\theta(x))\dot{R}_k\left(\frac{c_{n}\chi(x,y)+k_N(x,y)}{d(x,y)^{n-1}}\right)d\textrm{vol}_g(y)\]
where we understand that there is a smooth cut-off adapted to balls of radius $\rho$ and $\epsilon$. The cutt off function has non zero slope between $\frac{1}{2}\rho$ and $\rho$ and between $\epsilon$ and $\frac{3}{2}\epsilon$. We suppose $\epsilon>\rho$, otherwise this integral does not appear and might be assumed to be zero. Since we got at our disposal the logarithmic modulus of continuity, we have that
\[ |\theta(y)-\theta(x)| \leq \delta \ \text{for $d(x,y)< \epsilon$} \]
where $\delta>0$ should be taken in a proper way later on (which is in fact possible by choosing $\epsilon$ sufficiently small). Hence integrating by parts and using the aforementioned property we get:
\[|I_{\textrm{med}}(x)||\nabla_{g}\theta|^{2}\leq C_{2} \frac{\delta}{\rho}|\nabla_{g}\theta|^{2}+O\left(\norm{\nabla_{g}\theta}^{2}_{L^{\infty}}\right)\]
with $C_{2}=C_{2}(\epsilon,M,\chi,\norm{\theta_{0}}_{L^{\infty}})$.
To make sure that the term $|I_{\textrm{med}}(x)||\nabla_{g}\theta|^{2}$ does not exceed the cubic term of (\ref{test}), we choose 
\[ \delta=\frac{C}{\norm{\theta}^{2}_{L^{\infty}}} \]
with a proper constant $C$. The outer piece can be bounded easily throughout integration by parts by
\[ |I_{\textrm{out}}(x)||\nabla_{g}\theta|^{2} \leq C_{3}|\nabla_{g}\theta|^{2} + O\left(\norm{\nabla_{g}\theta}^{2}_{L^{\infty}} \right) \]
where $C_{3}= C_{3}(\norm{\theta_{0}}_{L^{\infty}},\chi,\epsilon)$. \\ \\
This proves the claimed inequality (\ref{test}):
\[  \frac{1}{2}L(|\nabla_{g}\theta|^{2})(x)+\frac{|\nabla_g \theta(x)|^{3}}{c\norm{\theta}_{L^{\infty}}} \leq C |\nabla_{g}\theta|^{2}(x)+ O(\norm{\nabla_{g}\theta}^{2}_{L^{\infty}}).\]

Finally, we close the argument: let $K$ be bigger than $C\|\theta_0\|_{\infty}$ from the hypothesis of Proposition \ref{muesca}, $2\|\nabla\theta_0\|_{\infty}$ and a constant such that the cubic term absorbs the quadratic in the right hand side above. Notice $K$ can be taken independent of $\nu$.  This choice of $K$ is enough to show that $|\nabla_g\theta|^2\leq K$ a fact that will be proved by contradiction. 

Let $x(t)\in M$ be the point where $|\nabla_g\theta|^2(\cdot, t)$ attains its maximum (this is well defined due to continuity and compactness). Even though $t\mapsto x(t)$ is not necesarily continuous, $t\mapsto|\nabla_g\theta|^2(x(t),t)$ is. As a consequence
\[t_0=\inf\{t\in(0,\infty):|\nabla_g\theta|^2(x(t),t)\geq K\}\textrm{ is positive}\]
We want to prove that $t_0=\infty$. If not, there exists some finite $t_0>0$ for which, by continuity, $|\nabla_g\theta|^2(x(t_0),t_0)\geq K$. But by definition one also knows $|\nabla_g\theta|^2(x(t_0),t)<K$ for any $t<t_0$. This facts altogether imply
\[\frac{d}{dt}|\nabla_g\theta|^2(x(t_0),t_0)\geq 0,\]
 which contradicts the inequality above read at the maximum $(x(t_0),t_0)$.

\section{\textbf{Appendix}}\label{s:appendix}

Here we will deliberately provide a sketchy idea of the proof which, nevertheless, follows a well-known pattern base on energy estimates. We include this appendix for completeness. 

\subsection{Local existence in $H^{m}(M)$ for $m>\frac{3}{2}$}
The divergence free condition $\textrm{div}_{g} u=0$ gives
\[\frac{1}{2}\frac{d}{dt}\norm{\theta}^{2}_{L^{2}(M)}\leq 0\]
which implies the decay of the $L^2$-norm. Let us denote by $D^n$ some arbitrary $n$th order derivative. Computing the evolution of the $H^{3}$ norm: 
\[ \frac{1}{2}\frac{d}{dt}\norm{D^3\theta}^{2}_{L^{2}(M)} + \norm{\Lambda^{\frac{1}{2}}(D^{3}\theta)}_{L^{2}(M)}  = -\int_{M} D^{3}(u\cdot \nabla_{g}\theta)D^{3}\theta \  + \int_{M} [\Lambda,D^{3}]\theta \ D^{3}\theta\]
The right hand side contains the nonlinear terms which can be expressed basically as four different terms. We will handle each one separately: if all derivatives fall with $u$
\[\int_{M} |D^{3}u \ \nabla_{g}\theta \ D^{3}\theta|\leq \norm{D^{3}\theta}^{2}_{L^{2}(M)}\norm{\nabla_{g}\theta}_{L^{\infty}(M)}; \]
if just two of them do
\[\int_{M}|D^{2}u \ D \nabla_{g} \theta \ D^{3}\theta| \leq \norm{D^{2}\theta}^{2}_{L^{4}(M)}\norm{D^{3}\theta}_{L^{2}(M)}\leq \norm{D^{3}\theta}_{L^{2}(M)}\norm{D^{2}\theta}_{L^{2}(M)}\norm{D^{3}\theta}_{L^{2}(M)} \]
\[\int_{M}|Du \ D^{2}\nabla_{g}\theta \ D^{3}\theta| \leq \norm{Du}_{L^{\infty}(M)} \norm{D^{3}\theta}^{2}_{L^{2}(M)};\]
if just one does
\[\int_{M}|[D^{3},\nabla_{g}]\theta u|\leq \norm{[D^{3},\nabla_{g}]\theta}_{L^{2}(M)}\norm{u}_{L^{2}(M)} ;\]
The last term corresponds to the commutator which can be bounded as follows
\[ \int_{M} [\Lambda,D^{3}]\theta \ D^{3}\theta \leq \norm{[\Lambda,D^{3}]\theta}_{L^{2}(M)}\norm{D^{3}\theta}_{L^{2}(M)} \]
All the bounds above imply
\[\frac{d}{dt}\norm{D^{3}\theta}^{2}_{L^{2}(M)} \leq C \norm{D^{3}\theta}^{3}_{L^{2}(M)}\]
where $C=C(M,\|\theta_{0}\|_{L^2})$. This inequality implies local in time boundedness for $\norm{\theta}_{H^{3}(M)}$. One may proceed in a similar fashion for higher order Sobolev spaces.\\

\textsc{Remark:} a standard application of Kato-Ponce commutator estimates (which have been extended for compact manifolds and more general settings, \cite{TAY, KP}) improves the local existence result to $H^{s}$ for any $s>\frac{3}{2}$. Notice that during the time of existence $T$,  the equation implies that \[\int_{0}^{T}\int_{M} |\Lambda^{s+\frac{1}{2}}\theta|^{2} < \infty.\]
As a consequence we can replace our initial time by some $t_{0}$ as close as we want to the initial time so that $\theta(x,t_{0})\in H^{2+\epsilon}$ for some $\epsilon>0$. Sobolev embedding implies now that $\nabla \theta(\cdot,t_{0})\in L^{\infty}$ which is all we need.


\subsection{Conditional global existence in $H^m(M)$}

In this section we show how global control of $\norm{\nabla_g\theta}_{L^{\infty}}$ provides inmediately global well posedness (cf. Theorem \ref{global}) in higher order Sobolev spaces $H^m$ for $m\in\mathbb{N}$. We do employ the critical  disipassion and since the case $m=0,1$ are easier to handle directly we will provide a sketch with $m>1$.

Consider the estimate for the $H^{m}$ norm of the solution
\[
 \frac{1}{2}\frac{d}{dt}\norm{\theta}^{2}_{H^{m}(M)}+ \norm{\theta}^{2}_{H^{m+\frac{1}{2}}(M)}  \leq \left|\int_{M} (u\cdot \nabla_{g}\theta)(-\Delta)^{m}\theta \ d\textrm{vol}_{g}(x)\right|
\]
Assume without loss of generality,that $m>1$ is an integer. Let us examine the nonlinear term more carefully, which can bounded using a commutator estimate by
\[C\int u\cdot[D^m, \nabla_g]\theta D^m\theta d\textrm{vol}_g(x)\]
and
\[ C \sum_{j=1}^{m} \int_{M} |D^{j}u||D^{m-j+1}\theta||D^{m}\theta| \ d\textrm{vol}_{g} (x) \]
The first might be bounded taking into account the fact that $u\in L^3$ uniformly. Appropiate H\"older and Sobolev inequalities dispose of this case (see below). For the sake of simplicity let us estimate it first for $j=1$ indicating later how to proceed in the general case. The rest of terms can be bounded easily in a similar way. Using the H\"older inequality and the boundedness of the Riesz transform in $L^{p}(M)$ spaces for $1<p<\infty$, yield
\[\begin{split}
 \int_{M} |Du||D^{m}\theta||D^{m}\theta| \ d\textrm{vol}_{g}(x) &\leq C \norm{Du}_{L^{3}(M)}\norm{D^{m}\theta}^{2}_{L^{3}(M)}\\
 &\leq \norm{D\theta}_{L^{3}(M)}\norm{D^{m}\theta}^{2}_{L^{3}(M)} \\
 \end{split} \]
Again H\"older and fractional Sobolev inequality imply
\[\begin{split}
\norm{D^{m}\theta}^{2}_{L^{3}(M)} &\leq \left(\int_{M} |D^{m}\theta|^{2} d\textrm{vol}_{g}(x)\right)^{\frac{1}{3}} \left(\int_{M} |D^{m}\theta|^{4} \ d\textrm{vol}_{g}(x)\right)^{\frac{1}{3}}\\
&\leq \norm{D^{m}\theta}^{\frac{2}{3}}_{L^{2}(M)}\norm{D^{m}\theta}^{\frac{4}{3}}_{L^{4}(M)} \\
&\leq \norm{D^{m}\theta}^{\frac{2}{3}}_{L^{2}(M)}\norm{D^{m+\frac{1}{2}}\theta}^{\frac{4}{3}}_{L^{2}(M)}
\end{split}
\]
and Galiardo-Nirenberg's inequality yields
\[ \norm{D^{m}\theta}_{L^{2}(M)} \leq C \norm{D\theta}_{L^{2}(M)}^{\frac{1}{2m-1}}\norm{D^{m+1/2}\theta}^{1-\frac{1}{2m-1}}_{L^{2}(M)}.\]
Therefore, the nonlinear term is bounded by
\[\begin{split}
 \int_{M} |Du||D^{m}\theta||D^{m}\theta| \ d\textrm{vol}_{g}(x) &\leq C \norm{D\theta}_{L^{3}(M)}\norm{D\theta}^{\frac{2}{3(2m-1)}}_{L^{2}(M)}\norm{D^{m+\frac{1}{2}}\theta}^{\frac{4}{3}+(1-\frac{1}{2m-1})\frac{2}{3}}_{L^{2}(M)} \\
 &\leq  C \norm{D\theta}_{L^{3}(M)}\norm{D\theta}^{\frac{2}{3(2m-1)}}_{L^{2}(M)}\norm{D^{m+\frac{1}{2}}\theta}^{2-\frac{2}{3(2m-1)}}_{L^{2}(M)}.
  \end{split}
\]
Since we have $\norm{\nabla_{g}\theta}_{L^{\infty}(M)}$ uniformly bounded and $M$ is a compact manifold, the above gives
\[\frac{1}{2}\frac{d}{dt}\norm{\theta}^{2}_{H^{m}(M)} + \norm{\theta}^{2}_{H^{m+\frac{1}{2}}(M)}\leq C \norm{\theta}^{2-\frac{2}{3(2m-1)}}_{H^{m+\frac{1}{2}}(M)}.\]
Finally, one may use Young's inequality to subsume in the dissipation term the right hand side. As a consequence $\frac{d}{dt}\|\theta\|^2_{H^m}\leq C$ where the $C=C(m,M,\norm{\nabla_{g}\theta}_{L^{\infty}})$.  

For the general case, $j>1$, the same strategy would work taking into account the following Gagliardo-Nirenberg inequality
\[ \norm{D^{m}\theta}_{L^{4}(M)} \leq \norm{D^{m}\theta}^{\frac{1}{2}}_{L^{2}(M)}\norm{D^{m+1}\theta}^{\frac{1}{2}}_{L^{2}(M)} \]
Hence the only case where we can not apply directly this estimation is when we have two terms with all the derivatives, which just occurs if $j=1$. 

\subsection{Uniqueness for smooth enough initial data}

Any two weak solutions $\theta_1$ and $\theta_2$ coincide provided that one of them is smooth. We denote the corresponding velocity vectors by $u_1$ and $u_2$, respectively. Let us define $f=\theta_1-\theta_2$, the variation of its $L^2$ norm can be estimated as follows
\[\frac{1}{2}\frac{d}{dt}\|f\|_{L^2(M)}^2=\frac{1}{2}\int_Mu_1\cdot\nabla_g(f^2)d\textrm{vol}_g+\int_M(u_1-u_2)\cdot\nabla_g\theta_2 fd\textrm{vol}_g-\|f\|^2_{H^{1/2}(M)}\]
the first term in the right hand side vanishes due to the incompressibility of $u_1$, the third might be neglected while the second might be bounded by
\[\|\nabla_g\theta_2\|_{\infty}\|f\|^2_{L^2(M)}\]
using the boundedness of the Riesz transforms and the hypothesis of $\theta_2$ being smooth for all times. Gronwall's inequality yields $f=0$ and hence the uniqueness.

\section{\textbf{Acknowledgments}}

The authors are indebted to T. Pernas-Casta\~no who carefully draw the figures herein. The third named author is grateful to V. Vicol who pointed out some helpful references in early stages of this work.

This work has been partially supported by ICMAT Severo Ochoa project SEV-2015-0554 and the MTM2011-2281 project of the MCINN (Spain).

\end{document}